# Quantifying USA tariffs effect: machine learning, entropy and fractal insights into the stock markets

**Mridul Patel**[1]

**Abstract**

This study presents a multiscale econometric framework to evaluate the impact of the USA tariff announcement of 2 April 2025 on the S&P 500 (USA) and S&P/ASX 200 (AUS) stock indices. We employ normalized permutation entropy (PE) to characterise the evolution of ordinal-pattern complexity and compute fractal-dimension estimates (Higuchi, Katz, Sevcik) to assess geometric scaling behaviour across different time-windows. Post-event PE remains uniformly high across windows, with values ranging from 0.666 to 0.913 for the USA and 0.690 to 0.923 for AUS, implying stable distributional entropy and no significant alteration of underlying symbolic dynamics. Fractal dimension estimates exhibit small but systematic changes: the Higuchi dimension increases (USA: 1.507 → 1.561; AUS: 1.524 → 1.533), indicating a marginal rise in high-frequency roughness, while Katz and Sevcik dimensions decline (USA Sevcik: 1.337 → 1.229), consistent with smoother medium-scale structure. Across all measures, the shock does not generate a statistically meaningful structural break. Additionally, machine-learning models (kNN, SVR, Random Forest, XGBoost, Neural Network) demonstrate strong cross-market predictability, with ensemble methods achieving out-of-sample $R^2 > 0.98$. Overall, the results suggest that both markets absorbed the tariff shock rapidly, exhibiting stable multiscale dynamics despite heightened geopolitical uncertainty. A model-agnostic XAI framework combining SHAP attribution and permutation-based diagnostics is used to isolate robust and independent information content.

✉ Mridul Patel
mridul.patel@student.rmit.edu.au; mridul.patel9999@gmail.com

[1] School of Mathematical Sciences, RMIT University, Melbourne, Australia

## 1 Introduction

Tariffs imposed by USA President Donald Trump beginning in early 2025 are projected to significantly reshape global trade dynamics and introduce substantial uncertainty into financial markets. These policies, which include a baseline 10% tariff on all USA imports and higher prices on goods from countries such as Australia, India, China, Canada, and members of the European Union, aim to address trade imbalances and strengthen domestic manufacturing. However, economic projections indicate that these tariffs may lead to increased costs for businesses and consumers, thereby reducing corporate profitability and weakening consumer demand.

The stock market is expected to experience heightened volatility amid the uncertainty surrounding tariff implementations and retaliatory measures from key trading partners. On April 2, 2025, the announcement of tariffs triggered a sharp market sell-off, with the S&P 500 projected to lose approximately 10% of its value over two days, resulting in an estimated $6.6 trillion loss in global equity. Volatility in stock markets serves as a key indicator for evaluating fluctuations in asset prices. As global financial systems grow more complex, accurately forecasting such volatility has become a focal point for both investors and academic researchers Ayitey Junior et al. (2023).

Traditional econometric models such as the Autoregressive Conditional Heteroskedasticity (ARCH) and (GARCH) Generalized Autoregressive Conditional Heteroskedasticity frameworks have long been used to model financial volatility (Engle 1982; Bollerslev 1986). These models capture time-dependent patterns, including volatility clustering, to a certain degree. However, they often struggle with nonlinear relationships and high-dimensional datasets. Moreover, their adaptability to fast-evolving market conditions and large-scale data remains limited, necessitating more flexible approaches for robust volatility forecasting (Andersen et al. 2006; Tsay 2010; Jin and Xu 2025).

To overcome these limitations, recent studies have increasingly embraced machine learning and deep learning techniques, which provide greater flexibility for modeling complex and nonlinear interactions within financial data. These methods can uncover latent structures and offer improved predictive performance in highly volatile and data-rich environments (Barua et al. 2024).

Forecasting stock market behavior has expanded beyond merely predicting price movements. It now encompasses the impact of various factors, including investor sentiment, macroeconomic indicators, and policy interventions. The multidimensional and dynamic nature of financial data substantially increases the complexity of this task Tsay (2010). While models like ARCH and GARCH continue to offer tools for quantifying volatility, their effectiveness diminishes when applied to large datasets involving numerous interacting variables (Hansen and Lunde 2005; Brooks 2019). Advanced machine learning methods provide a compelling alternative, capable of capturing intricate patterns and nonlinear dynamics with higher efficacy Gu et al. (2020).

These advantages become especially salient during periods of market turbulence, when sudden shifts and nonlinear relationships dominate. In such environments, conventional statistical models frequently fail to capture the nuanced behavior of markets, resulting in reduced predictive accuracy. Consequently, researchers and practitioners have increasingly turned to adaptive and robust machine learning methods that can accommodate complexity and provide reliable forecasts under uncertain conditions.

Beyond its role in forecasting price trends, volatility modeling also informs key decisions in the pricing, risk management, and portfolio optimization of financial instruments. Enhanced forecasts therefore empower investors to develop better trading strategies and enable financial institutions to optimize capital allocation. Ultimately, more accurate volatility prediction reduces exposure to market risk and enhances return potential (Somkunwar et al. 2024; Xu and Zhang 2023a).

The k-Nearest Neighbors (kNN) algorithm has shown considerable potential in financial forecasting, predicting stock movements by comparing new data with historical patterns derived from features such as lagged returns and trading volume Cover and Hart (1967). Its non-parametric design enables flexible modeling of nonlinear relationships without distributional assumptions Altman (1992), though its performance declines in high-dimensional settings and requires careful preprocessing to mitigate noise in financial data Nti et al. (2020).

Random Forest Regression employs an ensemble of decision trees to predict market movements using a wide range of macroeconomic and technical variables (Breiman 2001). It is robust to noise and provides interpretable insights through feature importance analysis (Gu et al. 2020), yet its computational intensity and tendency to overfit under high volatility remain limitations (Sonkavde et al. 2023). Careful hyperparameter tuning and cross-validation are therefore essential for reliable financial forecasting.

Recent studies highlight the growing role of machine learning in financial forecasting. Gu et al. (2020) demonstrate that models such as Random Forests and neural networks outperform traditional linear approaches in asset pricing, emphasizing the importance of feature selection. Nti et al. (2020) review algorithms like kNN and SVR, noting their ability to capture nonlinear patterns but also their sensitivity to noise. Building on this, Hu et al. (2023) employ SVR with sentiment and technical indicators for Chinese A-share prediction, while Sonkavde et al. (2023) find Random Forests effective for volatility forecasting, albeit prone to overfitting in turbulent markets.

## 1.1 Existing literature on USA tariffs

A growing strand of empirical research investigates how tariff policies and trade uncertainty transmit to financial and macroeconomic outcomes. The effects of USA tariffs and trade policy uncertainty on major stock indices using a structural vector autoregression (SVAR) model with monthly data Yilmazkuday (2025). Their results indicate that one standard-deviation shocks to tariffs or trade policy uncertainty cause significant long-run reductions in stock prices between 7.3 and 10.1% across the S&P 500, Nasdaq Composite, and Dow Jones indices, with tariffs and uncertainty jointly explaining up to 9.9% of forecast error variance.

Similarly, Bandyopadhyay et al. (2025) employs SVARs to evaluate the broader macroeconomic impact of tariffs on the USA economy, showing that persistent tariff increases depress real GDP growth, reduce stock market capitalization, and worsen the fiscal balance, while generating only transitory inflationary effects. Further asymmetric responses of stock returns to tariff hikes and cuts using a nonlinear autoregressive distributed-lag (NARDL) framework Gaies (2025). Study finds that tariff hikes trigger immediate declines in returns followed by partial recovery, whereas both hikes and cuts lower returns in the long run, confirming that trade policy shifts are perceived as lasting sources of uncertainty.

Recent studies have also explored the global financial repercussions of tariff shocks. The study in Benguria and Saffie (2025) analyzes the April 2025 USA tariff announcement, distinguishing between anticipated and unanticipated components. They report that the "round-up" component has a disproportionately large effect, with a one-percentage-point tariff increase reducing global stock prices by 0.23% on average and modestly depreciating the USA dollar. Fiscal transmission channels, showing that when tariff revenues are used to subsidize investment or reduce labor taxes, short-term economic activity increases, and the welfare-maximizing unilateral tariff approaches 18% Alessandria et al. (2025).

### 1.2 Research gaps

While a substantial body of research has examined the financial and macroeconomic implications of tariff policies, most existing studies have relied on econometric approaches models to trace causal and asymmetric responses in prices, returns, and output. These analyses provide valuable insights into aggregate dynamics but often overlook the intrinsic complexity of market reactions and the evolving structure of financial time series following trade shocks.

To date, no work has jointly applied complexity-based measures, such as permutation entropy and fractal dimensions, alongside advanced machine learning regressors to evaluate short- and medium-horizon responses of equity markets to tariff announcements.

This study aims to fill existing gaps in the literature by combining market complexity analysis with nonlinear regression techniques. This combination allows for a better understanding of both structural irregularities and predictive behaviors in response to announcements regarding U.S. tariffs. By extending beyond traditional frameworks of volatility and causal inference, this approach provides a fresh perspective on how uncertainty in trade policy affects financial systems. The current literature highlights several significant gaps that this study seeks to address.

- Prior work on tariff announcements largely focuses on event-study returns or volatility measures, with limited attention to the structural complexity of price dynamics.
- No existing study has examined tariff-related market responses using permutation entropy or fractal-dimension estimators, despite their relevance for capturing nonlinear and multiscale behaviour.
- The combination of complexity measures with machine-learning regression models has not been explored in the context of trade-policy shocks.
- Existing studies do not compare pre- and post-event complexity to assess whether policy shocks alter multiscale market structure.
- Cross-market analysis using multiscale complexity tools specifically comparing the S&P 500 and S&P/ASX 200 remains unexplored in the literature.

### 1.3 Aims and objectives

This study aims to examine how major macroeconomic and geopolitical developments affect the behaviour of the S&P 500 and S&P/ASX 200 indices over the period from Nov-2020 to Nov-2025. The analysis focuses on the introduction of new US tariff measures on April 2, 2025, and evaluates whether these events generated measurable shifts in market dynamics. The objectives of the study are defined as follows:

1. To assemble a unified dataset of daily prices for the USA and AUS equity markets and characterise their statistical properties.
2. To quantify market complexity using permutation entropy and evaluate structural shifts through the change in entropy ($\Delta$PE) and corresponding standardized effect sizes.
3. To compute fractal dimensions using Higuchi's(1988), Katz's (1988), and Sevcik's (2010) estimators and assess how short-, medium-, and long-horizon geometric complexity evolved before and after the tariff announcement.
4. To model cross-market relationships using machine-learning methods and assess predictive accuracy across a range of alternative specifications.
5. To compare pre- and post-event behaviour in both indices and determine the extent to tariff announcement altered volatility patterns and multiscale market dynamics.

This study employs multiple regression models XGBoost, regression, Support Vector Regression (SVR), Random Forest Regression, k-Nearest Neighbors (kNN) regression, and Neural Network (MLP) to capture both linear and nonlinear dynamics in stock market behavior. Moreover, complexity measures such as permutation entropy and fractal-dimension estimators provide powerful tools for quantifying irregularity and structural variability in financial time series. Permutation entropy captures the degree of disorder in temporal patterns by evaluating the distribution of ordinal sequences, offering a dynamic measure of informational complexity in market behavior. Complementarily, fractal-dimension estimators including Higuchi's (1988), Katz's (1988), and Sevcik's (2010) methods characterize the geometric roughness of time series by relating signal detail to scale. Higuchi's method evaluates local scaling properties directly from the data, Katz's approach estimates signal complexity based on trajectory length and amplitude, and Sevcik's normalized curve-length measure enables efficient comparison across datasets. Together, these techniques capture both the informational and geometric aspects of market dynamics, providing complementary perspectives on financial complexity.

The structure of the paper is organized as follows: Section 2 describes the dataset, outlines the key events timeline, and discusses the principal drivers of stock market behavior. Section 3 details the econometric and regression-based methodologies employed for modeling and forecasting, and introduces the complexity measures such as permutation entropy and fractal dimensions used to quantify the effects of tariff shocks. Section 4 presents and interprets the empirical results obtained from the machine learning models and complexity analyses. Finally, Section 5 summarizes the main findings, highlights their broader implications, and acknowledges the limitations of the present study while suggesting directions for future research.

## 2 Data

This section describes the data sources employed in the analysis. Historical weekly price data for the S&P/ASX 200 and S&P 500 indices were obtained from publicly accessible financial repositories, specifically the official databases of the Australian Securities Exchange ASX

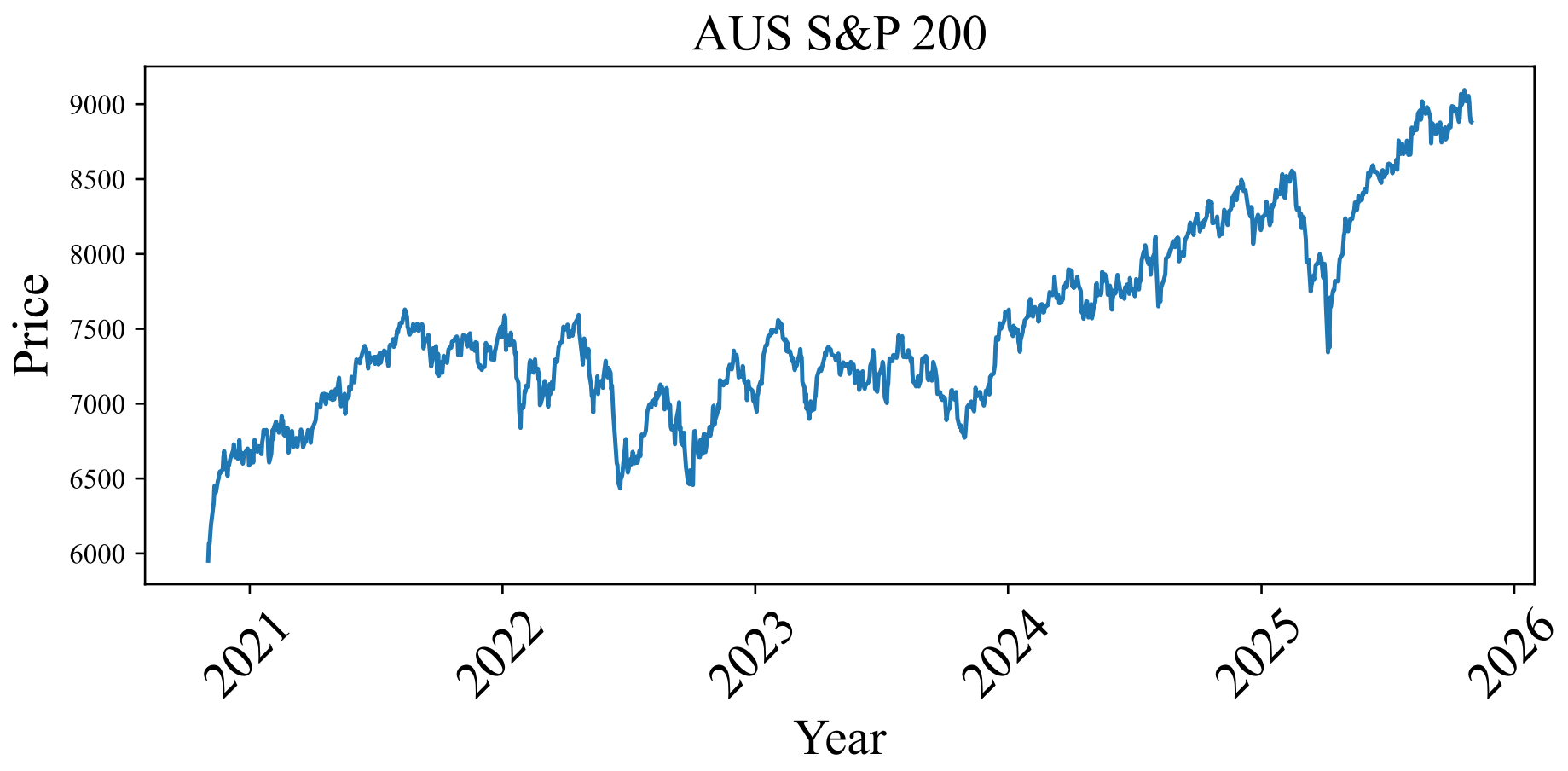

**Fig. 1** Daily closing prices of the AUS stock market over the same period. The price movement pattern reflects market behaviour that can be compared with the USA trends

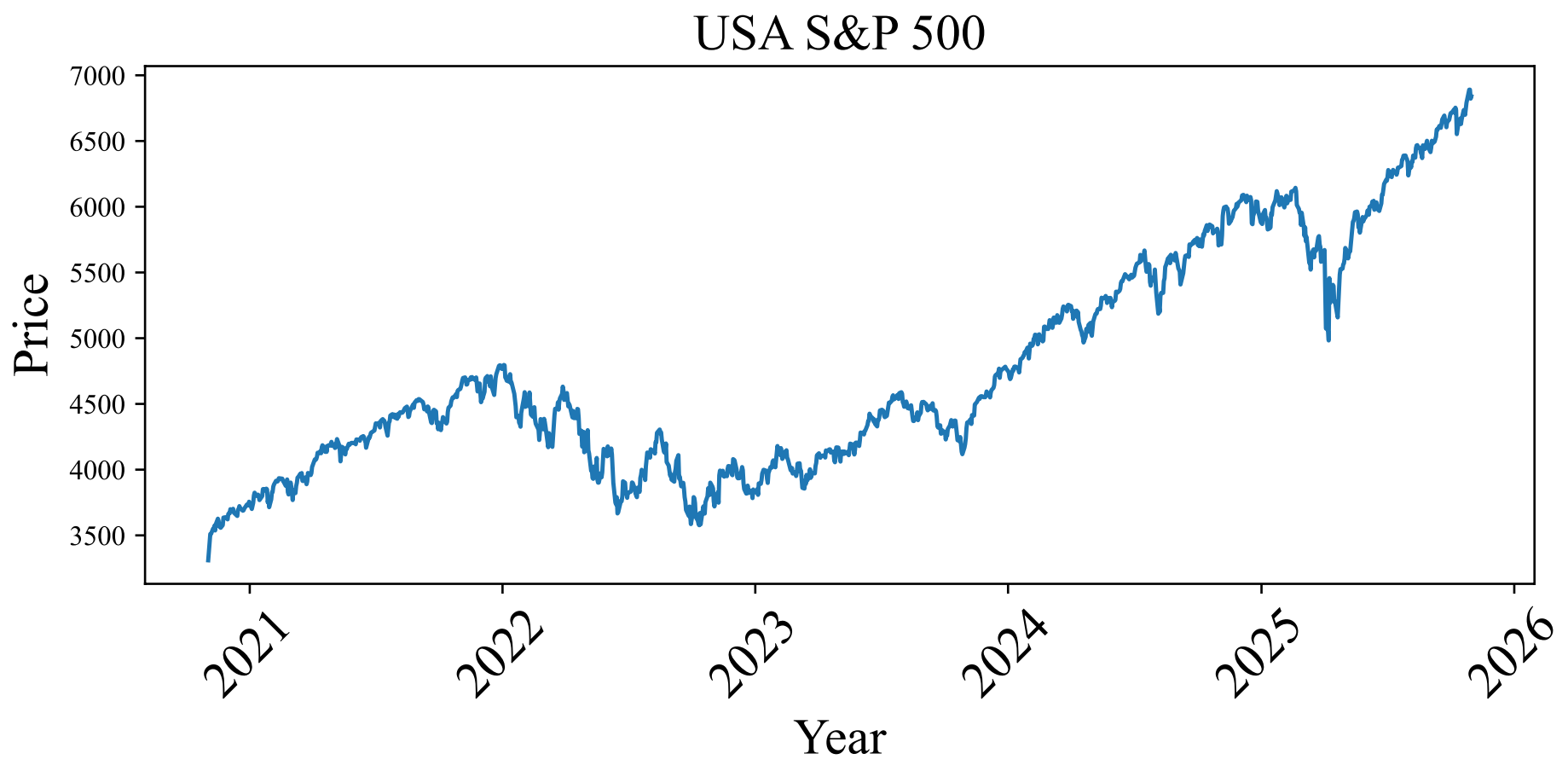

**Fig. 2** Daily closing prices of the USA stock market over the observed period. The plot captures short-term fluctuations and long-term trends based on historical data

Limited (2025) and S&P Dow Jones Indices S&P Dow Jones Indices (2025). The datasets compile the aggregated daily performance of the Australian stock market (S&P/ASX 200) and the USA stock market (S&P 500) across the period spanning November 1, 2020, to November 1, 2025. The daily closing prices of (S&P/ASX 200) and (S&P 500) are illustrated in Figs. 1 and 2, respectively.

The full sample spans the period from November 1, 2020, to November 1, 2025, a horizon that includes several key economic and geopolitical developments. These years capture prolonged pandemic-related uncertainty, the global recovery phase, shifts in monetary policy across major central banks, and the announcement and subsequent implementation of a broad tariff package in the USA, which came into effect on April 2, 2025. These events are reflected in the fluctuations of daily highs, lows, and closing prices in both datasets, enabling a clear comparison of pre- and post-event market behaviour during this period.

## 2.1 Key events timeline (2020–2025)

The dataset spans the period from November 1, 2020, to November 1, 2025, encompassing several major economic, political, and geopolitical developments that shaped global equity markets. This timeframe includes significant fluctuations affecting both the S&P 500 (USA) and the S&P/ASX 200 (AUS). The key events within this window are summarised below.

- *Late 2020–2021: Global pandemic conditions.* Continued COVID-19 disruptions, intermittent lockdowns, and fiscal stimulus measures shaped market volatility.
- *2021–2022: Monetary tightening cycle.* Beginning in late 2021, major central banks (Federal Reserve, Reserve Bank of Australia) initiated successive interest-rate increases in response to rising inflation, generating elevated uncertainty in both markets.
- *2023: Recession fears and inflation persistence.* Slower global growth and persistent inflationary pressures contributed to heightened risk aversion and episodic corrections in equity indices.
- *2024: Pre-election uncertainty in the USA.* The build-up to the USA presidential election produced short-lived volatility spikes in the S&P 500.
- *January 2025: Inauguration of Donald Trump's second term.* The start of the new administration was accompanied by shifts in policy expectations, particularly in areas related to trade, energy, and manufacturing.
- *April 2, 2025: Announcement and implementation of new USA tariff measures.* The introduction of broad tariffs on selected trading partners generated immediate reactions in global markets, reflected in short-horizon fluctuations in both the USA and AUS stock indices.
- *Mid–2025: Adjustment period post-tariffs.* Following the initial shock, markets exhibited a period of stabilisation, with prices adjusting to revised trade and supply-chain expectations (Fig. 3).

## 2.2 Key drivers of stock markets

Stock markets respond to a combination of macroeconomic conditions, financial fundamentals, and geopolitical developments that shape investor expectations and risk assessments. Core economic indicators such as GDP growth, inflation, employment, and interest rates signal the direction of economic activity and influence discounting and capital allocation decisions. Financial factors, including corporate earnings, sectoral performance, exchange rates, and commodity prices, further determine the outlook for firm profitability and investment flows.

Short-run market dynamics are also sensitive to shifts in liquidity, trading volume, and broader market sentiment, while geopolitical events ranging from elections and policy reforms to trade measures and international tensions can trigger rapid adjustments in prices through changes in perceived uncertainty. A graphical illustration of the key drivers of the stock market is presented in Fig. 4.

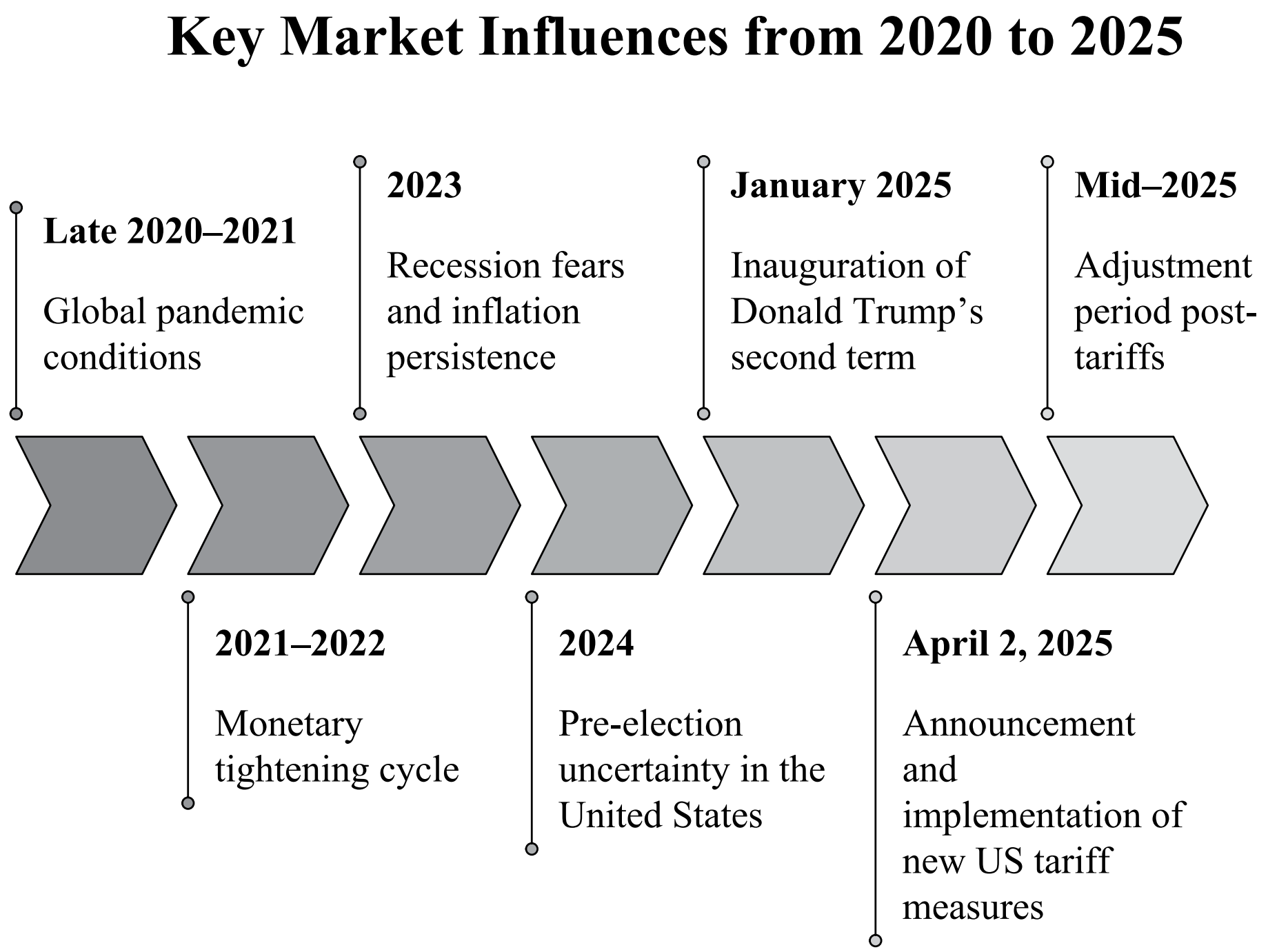

**Fig. 3** A timeline of major macroeconomic and geopolitical developments from 2020 to 2025, including pandemic disruptions, monetary-policy tightening, recession concerns, the 2024 USA election cycle, and the introduction of new USA tariff measures on April 2, 2025

## 3 Methodology

This section combines exploratory data analysis, statistical summaries, and machine-learning methods to support the empirical results. The EDA provides a clear view of the underlying structure of the data and highlights the main distributional features. Permutation entropy serves as a noise–tolerant indicator of time-series complexity by tracking the diversity of ordinal patterns. We compute both the change in the permutation entropy ($\Delta$PE) and its standardized effect size to evaluate the market sensitivity around key events. To examine geometric scaling across time horizons, we estimate fractal dimensions using the Higuchi, Katz, and Sevcik approaches. These measures show that post-event permutation entropy remains consistently elevated across all windows, signalling a sustained increase in dynamical complexity. The statistical analysis summarises descriptive characteristics and inter–variable relationships, while the machine-learning component outlines the predictive models employed for forecasting.

### 3.1 Regression models for stock market prediction

This section reviews the key regression models employed for stock market volatility and price prediction. Each model provides distinct advantages in capturing patterns within

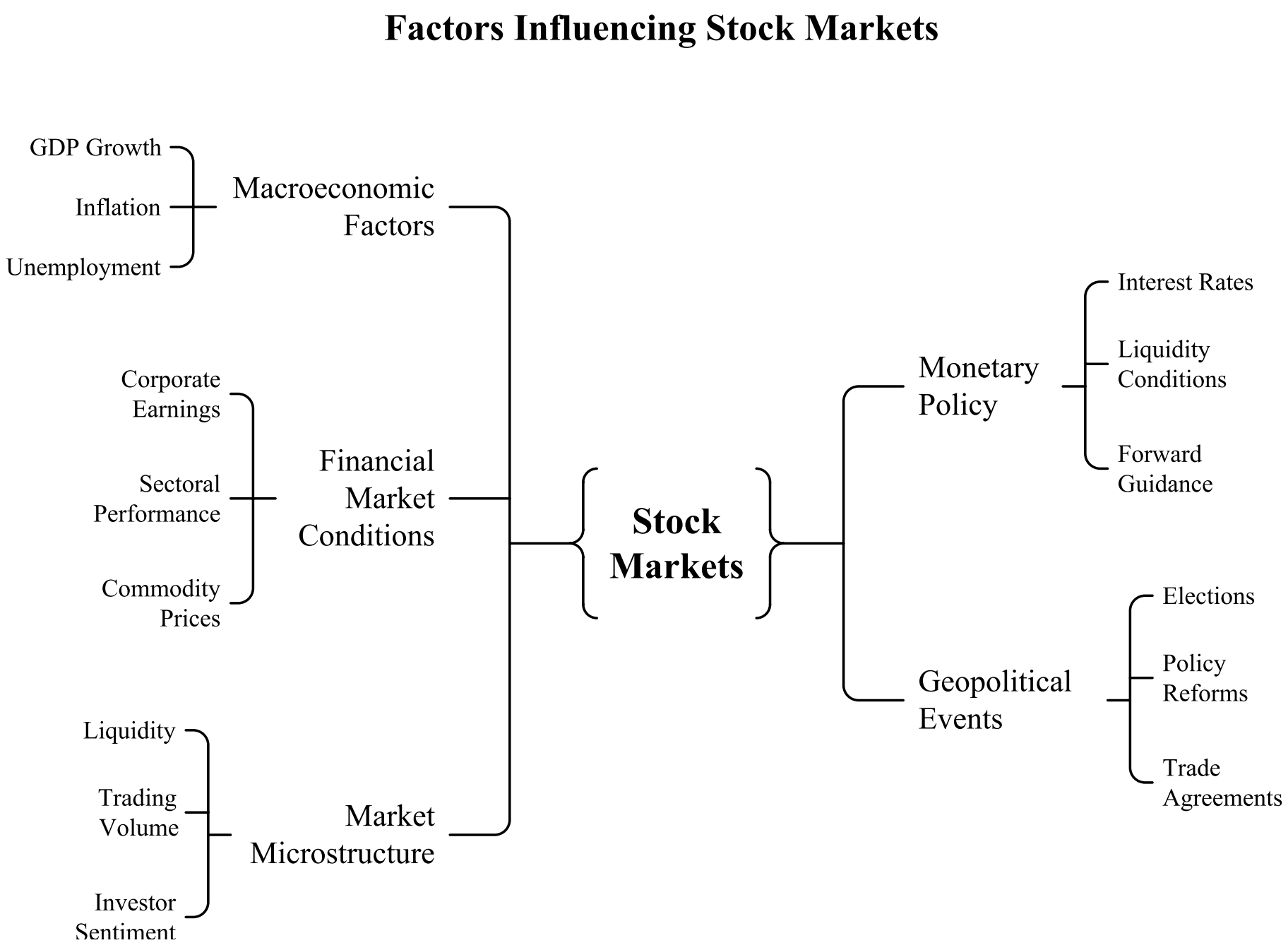

**Fig. 4** A conceptual illustration summarises the primary drivers of stock-market dynamics, including macroeconomic fundamentals, monetary-policy conditions, financial-market factors, and the influence of geopolitical events, policy shifts, and investor sentiment

historical financial data, ranging from simple linear relationships to complex non-linear dynamics.

### 3.2 Support vector machines (SVMs): linear and polynomial kernels

Support Vector Machines (SVMs) offer a robust framework for both linear and nonlinear regression tasks. When equipped with a linear kernel, SVMs are effective for high-dimensional financial data where relationships between variables are approximately linear. They provide strong regularization and generalization, making them suitable for structured prediction under noisy conditions. In contrast, SVMs with polynomial kernels extend the model's capacity to capture nonlinear patterns by projecting the input features into higher-dimensional spaces. This kernel-based flexibility allows the model to fit more complex financial behaviors such as regime shifts or nonlinear volatility patterns. Together, these variants of SVM provide a powerful toolkit for handling different levels of complexity in stock market data (Noble 2006; Basak et al. 2007).

### 3.3 Random forest regression

Random Forest Regression is an ensemble learning method that constructs multiple decision trees and aggregates their predictions to reduce variance and enhance robustness. It handles both linear and nonlinear relationships effectively, making it well-suited for finan-cial prediction tasks characterized by noise, non-stationarity, and complex interactions. Its ability to capture feature importance and resist overfitting contributes to its popularity in modeling stock returns, volatility, and macroeconomic indicators Liaw and Wiener (2002).

### 3.4 k-nearest neighbors (kNN) regression

k-Nearest Neighbors (kNN) regression predicts the target value for a new instance by averaging the outcomes of its $k$ closest neighbors in the feature space. As a non-parametric method, kNN does not make assumptions about the underlying data distribution and can model both linear and nonlinear trends. In financial applications, kNN is particularly valuable for short-term forecasting due to its sensitivity to local data structure, though its performance can degrade with high dimensionality and noisy features Altman (1992).

### 3.5 XGBoost regression

XGBoost is a gradient-boosting framework that builds ensembles of decision trees to capture nonlinear relationships and interaction effects in the data Chen (2016). Its use of regularisation and efficient optimisation makes it robust to overfitting and well-suited for modelling event-driven fluctuations in financial time series.

### 3.6 Neural network

The feed forward Neural Network is a flexible non-linear model that maps inputs to outputs through layers of interconnected neurons Gardner (1988). While it can learn rich, non-linear relationships, its performance depends heavily on the amount and quality of available data and can be sensitive to over-fitting, particularly in volatile market environments. For applications of neural networks reader may see (Xu and Zhang 2023c, b, a)

To better understand these volatility regimes and to complement the predictive insights from the ML regression models, we also examine the underlying complexity of the time series through permutation entropy. This metric provides a compact and noise-robust quantification of the diversity of ordinal patterns present in the data. Higher entropy values reflect greater irregularity and unpredictability, whereas lower values indicate more structured or deterministic behaviour.

In practice, shifts in entropy around major economic or geopolitical events help characterise how strongly market dynamics react to new information. To capture these changes, we compute both the variation in permutation entropy ($\Delta$PE) and its standardised effect size, defined as the deviation from the pre-event mean scaled by the corresponding pre-event standard deviation. These measures offer a complementary perspective to the ML forecasts by highlighting periods where market structure becomes more complex and potentially harder for nonlinear models to learn.

### 3.7 Computation of entropy

To quantify the time-varying complexity of the price series, we compute permutation entropy, a model-free measure that captures the diversity of ordinal patterns within a sliding window. Unlike traditional volatility metrics, which track the magnitude of price changes, permutation entropy reflects the underlying temporal structure of fluctuations and is therefore well suited for detecting subtle shifts in market organisation. By evaluating entropy across multiple window lengths, we obtain a multi-scale perspective on how informational complexity evolves around key economic events.

We quantify time-local complexity using permutation entropy (PE), following the Bandt–Pompe ordinal pattern approach Bandt (2005). Let $\{x_t\}_{t=1}^{T}$ be the financial time-series, where $t$ indexes time and $T$ is the total number of observations. For an embedding order $m = 3$ and a time delay $\tau = 1$, we construct delay vectors

$$\mathrm{x}_t^{(m,\tau)} = (x_{t-2}, x_{t-1}, x_t),$$

which represent sequences of $m$ consecutive observations spaced by $\tau$ time steps. Each delay vector is then mapped to its corresponding ordinal pattern (or permutation) $\pi_t \in S_m$, where $S_m$ denotes the set of all $m!$ possible orderings of $m$ elements.

Within a rolling window of length $w$ quarters ending at time $t$, we estimate the probability of observing a given pattern $\pi$ as

$$p_t(\pi; w) = \frac{N_t(\pi; w)}{\sum_{\pi' \in S_m} N_t(\pi'; w)},$$

where $N_t(\pi; w)$ is the number of times the pattern $\pi$ occurs in that window, and the denominator is the total number of patterns observed.

The normalized permutation entropy at time $t$ is then

$$\mathrm{PE}_{m,\tau}(t; w) = \frac{-\sum_{\pi \in S_m} p_t(\pi; w)\, \ln p_t(\pi; w)}{\ln(m!)} \in [0, 1],$$

where the numerator is the entropy of the pattern distribution, and the denominator $\ln(m!)$ ensures normalization so that the maximum possible entropy equals 1. A value close to 1 indicates high complexity and randomness in the sequence of patterns, while a value close to 0 indicates greater regularity and predictability.

### 3.8 Fractal dimensions

Fractals arise naturally in many physical and biological systems, including coastlines, clouds, and physiological signals, where structure is generated through recursive and non-linear mechanisms. Their prevalence in nature has motivated their use across diverse fields such as physiology Goldberger et al. (2002), finance Mandelbrot (1997), computer graphics Ebert et al. (2002), and geophysics Turcotte (1997). In financial economics, fractal methods provide a useful framework for characterising irregular time-series behaviour, examining volatility patterns Patel et al. (2025c), and quantifying market complexity (Verma and Kumar 2023; Patel et al. 2025a, b). They are equally relevant for fixed-income markets, where government bond yields encapsulate expectations about inflation, monetary policy, and macroeconomic conditions.

The fractal dimension offers a numerical measure of structural complexity and can take non-integer values Mandelbrot (1982), reflecting the degree of irregularity across scales.

In this study, we employ three widely used and computationally efficient estimators those due to Katz, Sevcik, and Higuchi each capturing a distinct aspect of fractal behaviour and performing well under rolling-window settings. More complex estimators (box-counting, wavelet-based methods, or correlation dimensions) are either computationally demanding or less stable on short time segments, making them less suitable for high-frequency analysis.

### 3.9 Katz's fractal dimension

Katz's method Katz (1988) relates the total curve length of a signal to its spatial extent. Let

$$P_i = (x_i, y_i), \qquad i = 1, \ldots, N,$$

denote sampled points, with Euclidean distance

$$d(P_i, P_j) = \sqrt{(x_i - x_j)^2 + (y_i - y_j)^2}.$$

The total curve length and diameter are

$$L = \sum_{i=1}^{N-1} d(P_i, P_{i+1}), \qquad d = \max_{1 \le i,j \le N} d(P_i, P_j).$$

Using the average step length $L/(N-1)$ for normalization, Katz's fractal dimension becomes

$$D = \frac{\log(N-1)}{\log(N-1) + \log\left(\frac{d}{L}\right)}.$$

### 3.10 Sevcik's fractal dimension

Sevcik's estimator Sevcik (2010) rescales the series to the unit square. For values $y_i$ sampled at equal intervals,

$$x_i^* = \frac{x_i}{x_{\max}}, \qquad y_i^* = \frac{y_i - y_{\min}}{y_{\max} - y_{\min}},$$

and the curve length in the normalized space is

$$L = \sum_{i=1}^{N-1} \sqrt{(x_{i+1}^* - x_i^*)^2 + (y_{i+1}^* - y_i^*)^2}.$$

The fractal dimension is then

$$D = 1 + \frac{\ln(L)}{\ln\big(2(N-1)\big)}.$$

### 3.11 Higuchi's fractal dimension

Higuchi's method Higuchi (1988) constructs subsampled sequences to quantify scaling behaviour. Given

$$Y = \{Y(1), \ldots, Y(N)\},$$

define for each $k = 1, \ldots, k_{\max}$ and $m = 1, \ldots, k$

$$Y_m^{(k)} = \left\{Y(m), Y(m+k), Y(m+2k), \ldots, Y\left(m + \lfloor \tfrac{N-m}{k} \rfloor k\right)\right\}.$$

The length of each subsequence is

$$L_m(k) = \frac{N-1}{\lfloor (N-m)/k \rfloor k^2} \sum_{i=1}^{\lfloor (N-m)/k \rfloor} |Y(m+ik) - Y(m+(i-1)k)|\,.$$

Averaging over $m$ yields $L(k)$, which scales as

$$L(k) \propto k^D,$$

so that the fractal dimension is estimated as the slope of the regression

$$\mathrm{HFD} = \frac{\log L(k)}{\log(1/k)}, \qquad k = 1, \ldots, k_{\max}.$$

We use $k_{\max} = 10$, a common choice that balances precision and computational cost.

Higuchi's estimator Higuchi (1988) captures the scaling behaviour of a signal and is particularly responsive to subtle local fluctuations. In contrast, Katz's method Katz (1988), which relies on the relationship between curve length and spatial extent, provides more stable estimates, especially when the data are short or noisy. Sevcik's approach Sevcik (2010), based on the length of the curve in a normalized coordinate system, yields dimension values that are directly comparable across series with different amplitudes or sampling ranges.

Together, these three estimators offer a practical balance between sensitivity, robustness, and computational efficiency, making them well-suited for rolling-window analysis of long financial time series.

### 3.12 Exploratory data analysis (EDA)

In this study, Exploratory Data Analysis (EDA) is conducted to investigate the structure and underlying characteristics of financial time series data. Key variables such as opening and closing prices, as well as volatility measures, are analyzed to gain insights into their tempo-ral dynamics. The primary aim is to uncover patterns, anomalies, or structural irregularities that may affect the performance and reliability of subsequent predictive modelling.

Time series plots are employed to visualize abrupt shifts and periods of heightened volatility, which are critical for understanding market dynamics. Furthermore, scatter plots and correlation heatmaps are used to explore interdependencies among variables, including the relationships between broad market indices and individual stock movements. During the EDA process, missing values and outliers are systematically identified. Appropriate imputation and filtering techniques are then applied to preserve data integrity. Based on the patterns revealed, relevant features such as lagged returns, technical indicators, and rolling statistical measures are extracted to enhance the feature set for machine learning models.

The descriptive statistics in Table 1 reveal clear differences between the two markets. The Australian index operates at a much higher level, with a mean of 7499.054 compared with 4754.875 for the USA, and this gap is consistent across the median and interquartile range.

However, volatility is noticeably greater in the USA, as its standard deviation (844.518) and empirical range (3580.690) both exceed those of AUS (617.719 and 3143.400, respectively), indicating wider fluctuations and stronger short-term variability. Both series are mildly right-skewed (USA: 0.694; AUS: 0.630), suggesting occasional upward movements, and both exhibit negative kurtosis (USA: −0.668; AUS: −0.251), pointing to distributions that are flatter and less heavy-tailed than a Gaussian benchmark.

The interquartile ranges reinforce these findings. In the USA series, the interquartile spread from Q1 (4114.555) to Q3 (5463.810) equals 1349.255, whereas the AUS spread from Q1 (7068.800) to Q3 (7850.700) is 781.900, almost 42% smaller. This demonstrates that middle-range variability is considerably higher in the USA market.

Taken together, these characteristics imply that while the Australian market is centred at a higher price level and remains comparatively stable, the USA market is more dispersed, more asymmetric, and more responsive to short-term disturbances patterns.

The relationship between the AUS and USA stock markets is examined using Pearson correlation on temporally aligned daily data. After preprocessing and merging by date, correlations are computed for key variables such as opening prices, closing prices, and returns. The strong positive correlations indicate synchronized market behavior, reflecting shared global economic influences, investor sentiment, and cross-market linkages.

**Table 1** Descriptive statistics and additional distributional measures for USA and AUS price series

| Statistic | USA | AUS |
|---|---|---|
| Count | 1256 | 1265 |
| Mean | 4754.875 | 7499.054 |
| Std. Dev | 844.518 | 617.719 |
| Min | 3310.200 | 5951.300 |
| 25% | 4114.555 | 7068.800 |
| Median | 4469.415 | 7354.400 |
| 75% | 5463.810 | 7850.700 |
| Max | 6890.890 | 9094.700 |
| Range | 3580.690 | 3143.400 |
| Variance | 713210.271 | 381576.393 |
| Skewness | 0.694 | 0.630 |
| Kurtosis | −0.668 | −0.251 |

These findings are visualized in a correlation heatmap, as shown in Fig. 5, where high correlation values underscore the degree of alignment between the two markets. The evidence of synchronized trends supports the inclusion of cross-market features in forecasting models, potentially enhancing their predictive accuracy and robustness.

To investigate the temporal dynamics of market behavior, a time series regression framework is employed, incorporating lagged variables and rolling statistical features. This approach captures short-term dependencies and localized trends in the data, thereby enhancing the ability to track and predict fluctuations in stock index movements over time.

### 3.13 Prediction of the data using machine learning

Machine learning has become an integral component of modern technology, enabling tasks such as pattern recognition and predictive analytics that are otherwise infeasible through traditional methods. Its strength lies in its ability to process high-dimensional and multivariate datasets, often within dynamic and uncertain environments. As the volume of available data increases, machine learning algorithms tend to improve both in predictive accuracy and computational efficiency. However, effective training requires access to large-scale datasets that are representative, unbiased, and of high quality. Furthermore, selecting the most appropriate algorithm typically involves manually evaluating multiple models process that is not only time-consuming but also prone to human error.

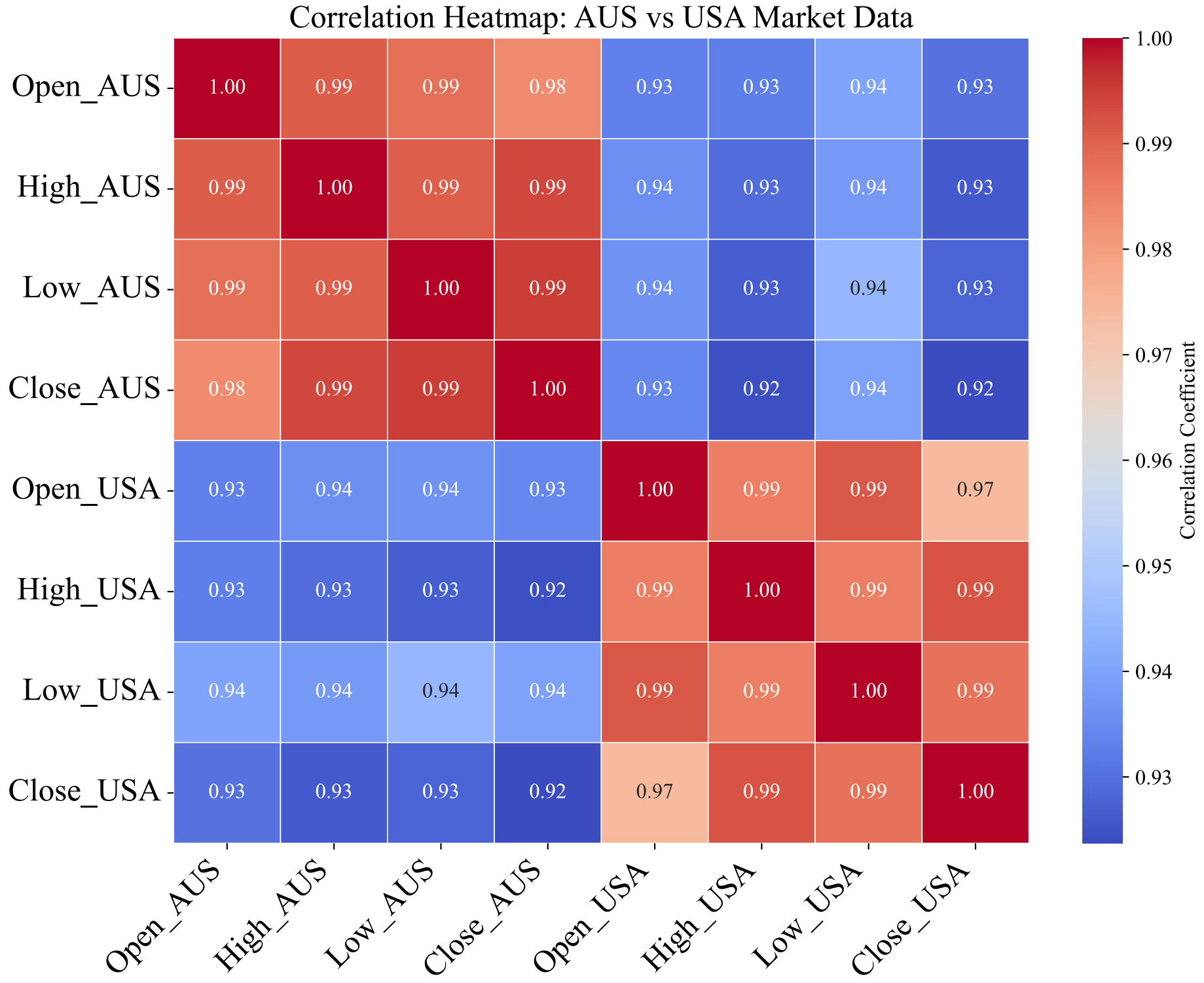

**Fig. 5** Correlation heatmap depicting the linear relationships between key stock market variables (Open, High, Low, Close) from the USA and AUS financial datasets

To investigate the predictive relationship between the Australian and USA stock markets, this study employs a range of machine learning models, including Support Vector Regression (SVR), K-Nearest Neighbors (KNN) Regression, Neural Network, XGBoost, and Random Forest Regression.

The models were trained using 80% of the dataset to predict stock market closing prices. After training, the remaining 20% of the data was used as the testing set to evaluate the predictive performance of each model. To assess the accuracy and quantify prediction errors, evaluation metrics including the coefficient of determination ($R^2$), mean squared error (MSE), and mean absolute error (MAE) were calculated for all forecasting models.

*Coefficient of determination* ($R^2$): The coefficient of determination, commonly known as $R^2$, measures the proportion of variance in the dependent variable that is predictable from the independent variables. It is a dimensionless metric bounded between 0 and 1, where values closer to 1 indicate better model fit. The formula for $R^2$ is given by:

$$R^2 = 1 - \frac{\sum_{i=1}^{N}(z_i - \hat{z}_i)^2}{\sum_{i=1}^{N}(z_i - \bar{z})^2}$$

where $z_i$ denotes the actual observed values, $\hat{z}_i$ represents the predicted values, and $\bar{z}$ is the mean of the observed values Draper and Smith (1998).

*Mean squared error (MSE):* The Mean Squared Error quantifies the average squared difference between the observed and predicted values, reflecting the variance of the residual errors. It is computed as:

$$\text{MSE} = \frac{1}{N}\sum_{i=1}^{N}(z_i - \hat{z}_i)^2$$

where $N$ is the number of observations. Lower MSE values indicate better predictive accuracy Hyndman and Athanasopoulos (2018).

*Mean absolute error (MAE):* The Mean Absolute Error represents the average of the absolute differences between observed and predicted values. It is calculated by:

$$\text{MAE} = \frac{1}{N}\sum_{i=1}^{N}|z_i - \hat{z}_i|$$

MAE provides an intuitive measure of average prediction error magnitude without considering direction Willmott (1985).

The relative error is calculated using the Mean Absolute Percentage Error (MAPE), which quantifies the average absolute difference between predicted and actual values as a percentage of the actual values. Formally, for a set of actual values $z_{\text{true}} = \{z_1, z_2, \ldots, z_n\}$ and predicted values $z_{\text{pred}} = \{\hat{z}_1, \hat{z}_2, \ldots, \hat{z}_n\}$, the Mean Absolute Percentage Error (MAPE) is defined as:

$$\text{MAPE} = \frac{100}{n}\sum_{i=1}^{n}\left|\frac{z_i - \hat{z}_i}{z_i}\right|$$

where $n$ denotes the number of observations. MAPE provides an interpretable metric of prediction accuracy in percentage terms, facilitating comparison across different models and data scales. It is important to note that MAPE may be sensitive to observations where $y_i$ is close to zero, which can lead to large relative errors. Despite this limitation, MAPE remains a widely used and intuitive measure for evaluating regression model performance in financial forecasting tasks.

To evaluate predictive performance, we compared the outputs of three regression models against actual stock market values over the analysis period. All models were trained on the same historical dataset to ensure consistency and fairness in evaluation. As shown in Fig. 6, the models demonstrate varying degrees of alignment with real market trends. The accompanying $R^2$ value and error metrics (MSE and MAE) provide a quantitative basis for assessing each model's forecasting accuracy and generalization capability, and is depicted in Fig. 7.

To investigate the relationship between the USA stock index (S&P 500) and the Australian stock index (S&P/ASX 200), a scatter diagram of daily closing prices is illustrated in Fig. 8. The results show a strong positive correlation, indicating that the market movements of these two indices are closely aligned throughout the observation period.

### 3.14 Economic interpretation of model results via explainable artificial intelligence (XAI)

While gradient-boosted tree models, such as XGBoost, demonstrate impressive predictive performance, their economic interpretability is often constrained by their black-box nature. High predictive accuracy alone does not indicate which variables provide meaningful economic insights, especially in the presence of strong multicollinearity and non-linear interactions. This limitation is particularly significant in financial and economic contexts, where understanding the explanatory factors and their policy implications is as crucial as achieving precise forecasts. To address this concern, we utilize Explainable Artificial Intel-ligence (XAI) techniques to differentiate true information content from redundant or correlated predictors.

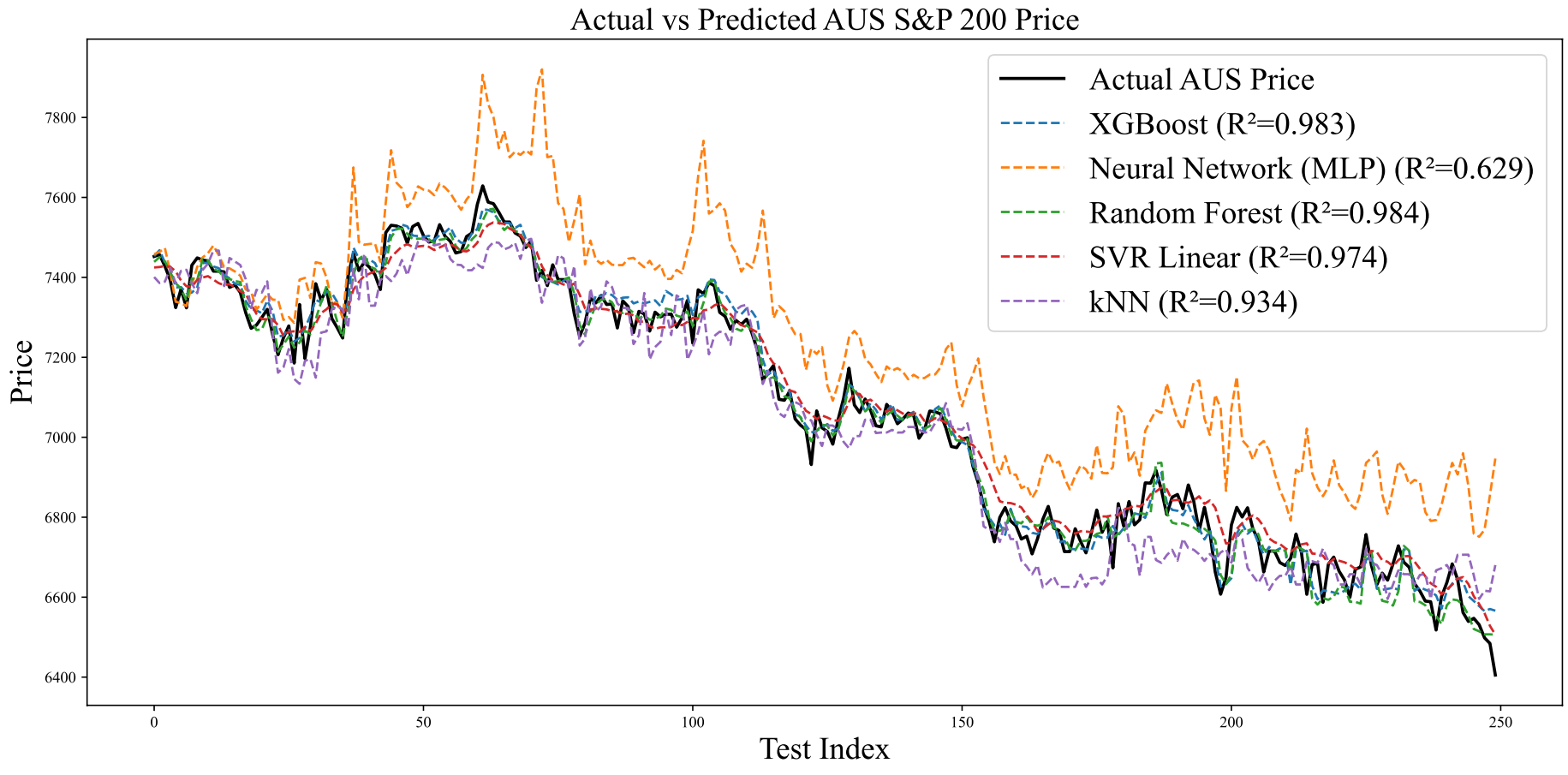

**Fig. 6** Comparison of actual and predicted stock market values using different models. The plot illustrates the performance of each model in tracking observed data over time, highlighting variations in predictive accuracy

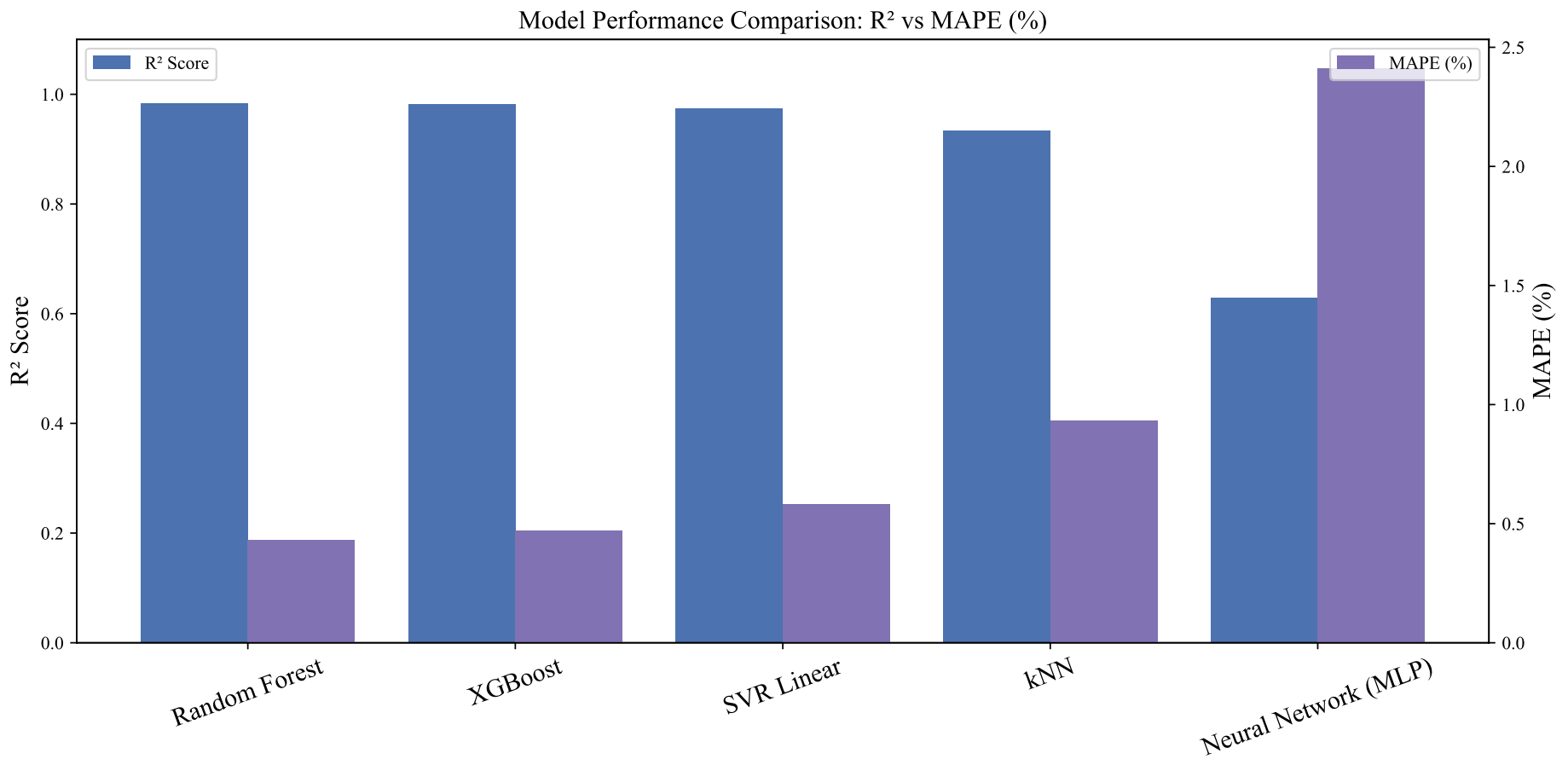

**Fig. 7** $R^2$ values and error metrics illustrating the predictive performance of different models

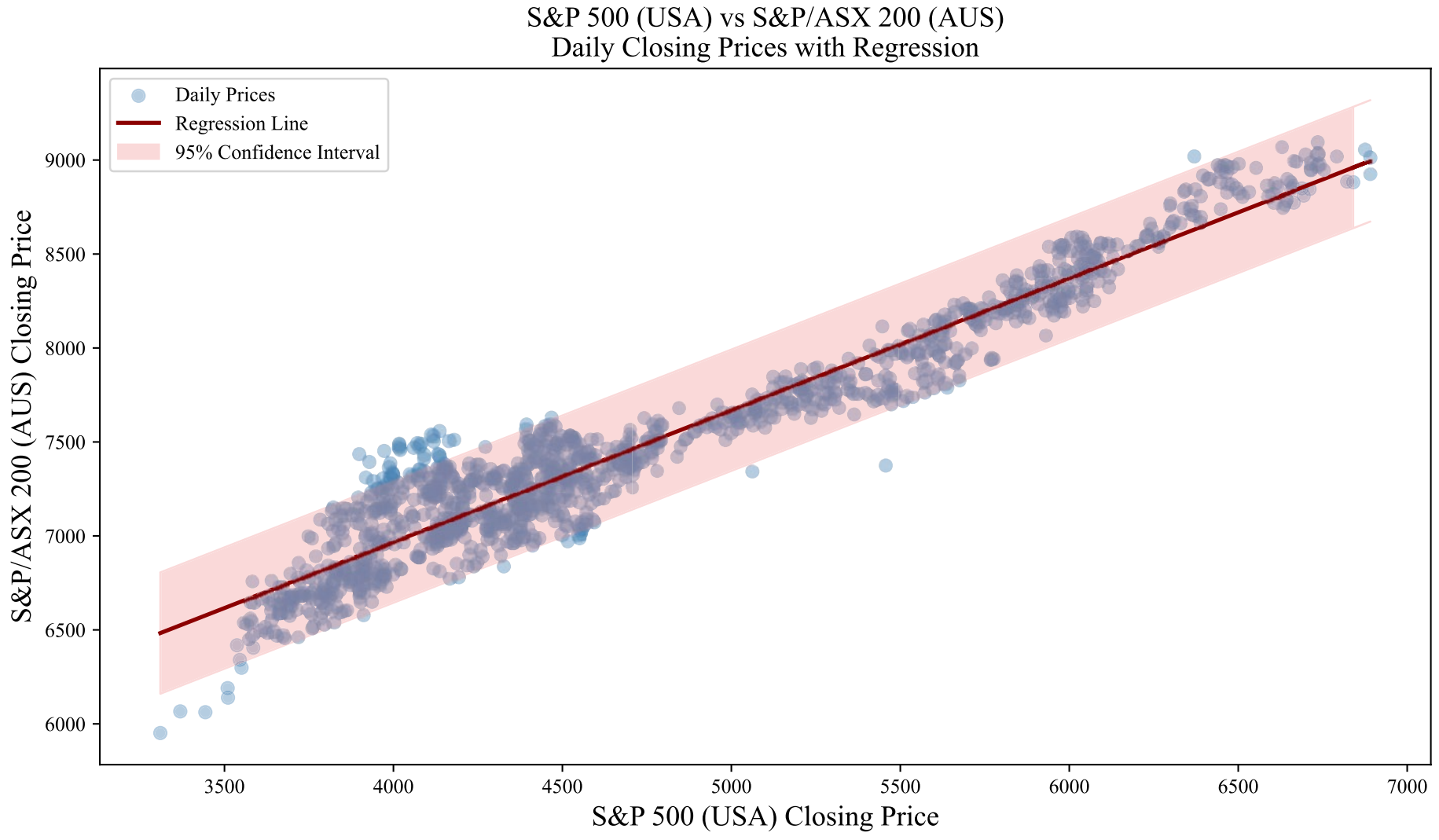

**Fig. 8** Scatter plot of daily closing prices for the S&P 500 (USA) and S&P/ASX 200 (AUS) indices with a fitted regression line and 95% confidence interval

Following the information-theoretic perspective of Berger Berger (2025), we move beyond raw feature importance and explicitly distinguish between (i) apparent importance driven by correlation, (ii) independent information content, and (iii) robustness of explanations. We adopt a model-agnostic SHAP framework applied to the prediction function of the trained XGBoost model, ensuring that interpretability does not rely on model-specific heuristics. In addition, we conduct redundancy diagnostics, conditional relevance analysis, and stability checks to obtain economically meaningful explanations.

### 3.14.1 Global feature importance

Table 2 reports the mean absolute SHAP values, which capture the average magnitude of each variable's contribution to the model's predictions. A clear pattern emerges: intraday price information dominates the explanatory structure of the model. In particular, the daily low and high prices exhibit markedly larger SHAP values than the opening price, percentage returns, or trading volume, indicating that the intraday price range carries the most informative signal for explaining variations in the predicted closing price.

The SHAP values are expressed in the units of the predicted variable, their absolute magnitudes are not directly comparable across markets and primarily reflect differences in index price levels rather than differences in the underlying explanatory structure. Conditional $\Delta R^2$ therefore provides a complementary measure by quantifying the independent information content of each predictor through the associated loss in out-of-sample fit when that feature is permuted. Variations in the magnitude of $\Delta R^2$ across markets are naturally expected and can be attributed to differences in volatility regimes and index construction, rather than to qualitative changes in the nature of the information driving price formation. The stability metric further reinforces this interpretation by assessing the robustness of SHAP-based importance across bootstrap resamples; higher stability values indicate that a variable's contribution is consistently identified and not driven by sample-specific fluctuations. Taken together, these diagnostics allow us to disentangle scale effects from genuine informational differences and to identify predictors that are both independently informative and reliably important across markets.

Table 2 shows that, for both AUS and the USA, the predictive structure of the model is overwhelmingly driven by intraday price extremes. Global SHAP values indicate that the daily high and low prices dominate feature importance, with all other variables contributing an order of magnitude less. Conditional relevance analysis sharpens this result: only the high and low prices generate substantial reductions in out-of-sample $R^2$ when permuted, while opening prices, percentage changes, and trading volume add virtually no independent information once the price range is known. Stability diagnostics further confirm that these findings are robust across resamples. Although the relative prominence of high versus low prices differs slightly between the two markets, the overall information structure is remarkably consistent, suggesting that intraday price range information is a structural driver of price formation rather than a market-specific artefact.

### 3.14.2 Redundancy and information overlap

Although the global SHAP rankings suggest that several intraday price variables are influential, a closer inspection reveals substantial redundancy among them. SHAP-based corre-lation diagnostics indicate extremely strong dependence between the High, Low, and Open prices, implying that these variables largely encode the same underlying information about the daily price range rather than contributing distinct explanatory signals.

**Table 2** Global SHAP importance, conditional relevance, and stability: AUS versus United States

| | AUS | | | USA | | |
|---|---|---|---|---|---|---|
| Feature | Mean \|SHAP\| | Cond. $\Delta R^2$ | Stability | Mean \|SHAP\| | Cond. $\Delta R^2$ | Stability |
| Low | 303.83 | 0.632 | 8.51 | 272.22 | 0.317 | 8.91 |
| High | 217.19 | 0.341 | 6.28 | 386.29 | 0.601 | 16.76 |
| Change (%) | 16.14 | 0.003 | 0.65 | 11.63 | 0.001 | 0.39 |
| Open | 15.34 | 0.003 | 0.36 | 44.01 | 0.008 | 2.14 |
| Volume | 1.65 | 0.000 | 0.07 | $2.6 \times 10^{-7}$ | 0.000 | $8.2 \times 10^{-7}$ |

#### 3.14.3 Conditional relevance and information content

To isolate genuine information content, we perform a conditional relevance analysis based on permutation-induced drops in out-of-sample $R^2$. The results, reported in Table 2, show that only the daily Low and High prices generate economically meaningful reductions in predictive performance when permuted. In contrast, opening prices, percentage changes, and trading volume have negligible conditional effects once the price range is accounted for. This finding demonstrates that the apparent importance of several covariates in the raw SHAP ranking is largely driven by correlation rather than independent information.

#### 3.14.4 Economic interpretation

Taken together, the XAI analysis reveals that the predictive success of the XGBoost model is primarily driven by information embedded in the intraday price range, rather than by trading volume or opening levels. While standard importance measures might suggest a broader set of influential variables. By explicitly separating importance, redundancy, and information content, the XAI framework ensures that the model's conclusions are economically meaningful and not artefacts of correlated predictors.

#### 3.14.5 Pre–post event interpretation around April 2, 2025

Table 3 reports mean absolute SHAP values computed over symmetric pre- and post-event windows for both AUS and the USA, together with the corresponding post-minus-pre changes. The target variable is the daily closing price in both markets.

The pre–post comparison indicates a clear shift in the intensity with which information is incorporated into prices following the April 2, 2025 event, while leaving the overall information hierarchy intact. In AUS, the post-event period is marked by a pronounced increase in the importance of the opening price, alongside a moderate rise in the role of trading volume, while the influence of intraday extremes declines relative to the pre-event window. This pattern suggests greater reliance on initial price signals immediately after the event, consistent with heightened uncertainty and reassessment of information at the market open.

A similar, though milder, adjustment is observed in the USA. The opening price again becomes more influential in the post-event period, the contribution of the high price increases, and the role of the low price remains largely unchanged. At the same time, percentage price changes lose relevance and trading volume continues to play a negligible role.

**Table 3** Pre–post event SHAP comparison around April 2, 2025

| | AUS | | | USA | | |
|---|---|---|---|---|---|---|
| Feature | Pre | Post | Change | Pre | Post | Change |
| Open | 54.35 | 131.99 | +77.64 | 94.97 | 156.66 | +61.69 |
| High | 118.89 | 97.23 | −21.66 | 44.44 | 64.88 | +20.44 |
| Low | 51.82 | 27.23 | −24.58 | 29.71 | 29.42 | −0.29 |
| Change (%) | 10.17 | 11.64 | +1.47 | 10.77 | 5.84 | −4.93 |
| Volume | 2.71 | 7.24 | +4.53 | $\sim 0$ | $\sim 0$ | $\sim 0$ |

Overall, the results indicate that the April 2, 2025, do not indicate a fundamental reordering of the information hierarchy underlying price formation. in either market. Instead, it temporarily reallocates attention toward opening-level information, reflecting an event-driven intensification of existing information channels rather than the emergence of new sources of economic information.

## 4 Results and discussion

We compare the predictive performance of six machine-learning models Random Forest, XGBoost, SVR, Linear SVR, kNN, and a Neural Network using an enhanced feature set that incorporates returns, rolling statistics, technical indicators, permutation entropy, and fractal-dimension measures. These features are designed to capture both conventional market patterns and deeper nonlinear structure. The resulting accuracy and error metrics for all models are summarised in Table 4.

The comparative performance of the five machine-learning models reveals important differences in their ability to capture the structure of daily stock-price dynamics. The Random Forest clearly emerges as the strongest performer, achieving the lowest prediction error (MAPE% $\approx 0.4323$) and the highest explanatory power ($R^2 = 0.984$). Its superior accuracy reflects the capacity of ensemble tree-based architectures to model nonlinear interactions, account for heterogeneous regimes in the data, and remain robust to noise and local irregularities. XGBoost performs nearly as well, with only slightly higher errors, confirming that boosted decision-tree methods effectively identify complex structure and successfully minimise residual variance through iterative refinement.

By contrast, the linear SVR, while achieving a respectably high $R^2$ of 0.974, displays larger absolute and percentage errors, signalling that a purely linear representation is inadequate for capturing the nonlinear components and volatility clusters typical of equity markets. The kNN model performs markedly worse, with its MAPE almost doubling relative to the SVR. This deterioration aligns with the well-known sensitivity of instance-based methods to noise, abrupt jumps, and changing local geometry features that are common in financial time series, especially around event-driven periods. The poorest performance is observed for the neural-network model ( MAPE% $\approx 2.41$, $R^2 \approx 0.6285$), where the combination of limited sample size, absence of temporal structure in the architecture, and potential overfitting lead to substantial prediction errors and unstable generalization.

Overall, the results underscore that ensemble tree models not only deliver the most reliable forecasts but also offer greater robustness across error metrics, highlighting their suitability for modelling short-horizon equity-price dynamics in environments subject to regime shifts, heteroskedasticity, and nonlinearity.

**Table 4** Model performance comparison for predicting daily stock prices using multiple machine-learning methods

| Model | MSE | MAE | $R^2$ | MAPE (%) |
|---|---|---|---|---|
| Random forest | 1619.85 | 30.15 | 0.9840 | 0.4323 |
| XGBoost | 1766.72 | 32.96 | 0.9826 | 0.4715 |
| SVR (linear) | 2667.11 | 40.82 | 0.9737 | 0.5827 |
| kNN | 6705.25 | 65.54 | 0.9338 | 0.9324 |
| Neural network | 37627.91 | 168.22 | 0.6285 | 2.4123 |

**Table 5** Change in permutation entropy (ΔPE) and standardized effect sizes for the USA (S&P 500) and AUS (S&P/ASX 200) following the tariff announcement on 2 April 2025

| Window | ΔPE (USA) | StdEffect (USA) | ΔPE (AUS) | StdEffect (AUS) |
|---|---|---|---|---|
| 8 | −0.0282 | −0.1724 | −0.0100 | −0.0604 |
| 12 | −0.0014 | −0.0117 | −0.0149 | −0.1199 |
| 16 | +0.0057 | +0.0597 | −0.0117 | −0.1250 |
| 20 | −0.0023 | −0.0306 | −0.0065 | −0.0893 |

**Table 6** Mean permutation entropy for the USA (S&P 500) and AUS (S&P/ASX 200) after the tariff announcement on 2 April 2025

| Window size | USA: mean PE | AUS: mean PE |
|---|---|---|
| 8 | 0.6658 | 0.6902 |
| 12 | 0.8078 | 0.7931 |
| 16 | 0.8656 | 0.8474 |
| 20 | 0.9129 | 0.9230 |

Although the machine-learning results highlight differences in predictive performance, they do not capture deeper structural features of market behaviour. To assess how the underlying price dynamics evolved around the tariff event, we complement the forecasting analysis with permutation entropy, a model-free measure that quantifies changes in temporal complexity and informational irregularity.

We evaluate $\mathrm{PE}_{3,1}(t;w)$ for $w \in \{8, 12, 16, 20\}$ using a one-quarter step, plotting each value at the midpoint of its corresponding window. As a result, visual responses to events appear with an alignment lag of approximately $w/2$ quarters (Table 5).

The estimated changes in permutation entropy reveal only limited adjustments in market complexity following the tariff announcement on 2 April 2025. Across all window lengths, the United States displays very small shifts (from $-0.0282$ to $+0.0057$), with standardized effect sizes close to zero, indicating that the event did not induce a substantive alteration in short-run informational irregularity. The Australian market shows uniformly negative changes in entropy (between $-0.0065$ and $-0.0149$), again with small standardized effects, suggesting a slight but minor reduction in unpredictability. Overall, the results point to weak and non-disruptive adjustments in both indices, consistent with the absence of a pronounced structural break in entropy-based complexity.

To characterise the level of informational complexity in both equity markets following the tariff announcement on 2 April 2025, we compute permutation entropy across multiple rolling-window lengths. This approach provides a multi-scale view of how the diversity of ordinal patterns evolves from short- to medium-horizon dynamics. The results, summarised in Table 6, show consistently high entropy values in both the S&P 500 and S&P/ASX 200, indicating that post-event price movements remained highly irregular and structurally comparable across markets.

The permutation-entropy estimates for the period after 2 April 2025 indicate consistently high levels of complexity in both equity markets. For the shortest window length ($w = 8$), the USA and Australia exhibit moderately elevated entropy values (0.666 and 0.690), pointing to irregular short-run price movements. As the window expands, entropy rises steadily, reaching 0.913 for the S&P 500 and 0.923 for the S&P/ASX 200 at $w = 20$, reflecting the increased mixing of ordinal patterns over longer horizons. The close alignment of the two markets across all window lengths suggests that, in the aftermath of the tariff announce-ment, both indices retained a similar degree of randomness and showed no evidence of structural divergence in their short- or medium-horizon informational complexity.

**Table 7** Fractal dimensions of the USA (S&P 500) and AUS (S&P/ASX 200) indices before and after the tariff announcement on 2 April 2025

| Market | FD Type | Before | After |
|---|---|---|---|
| USA | Higuchi | 1.5072 | 1.5612 |
| | Katz | 1.6187 | 1.5712 |
| | Sevcik | 1.3373 | 1.2290 |
| AUS | Higuchi | 1.5238 | 1.5329 |
| | Katz | 1.6408 | 1.5863 |
| | Sevcik | 1.3822 | 1.2475 |

**Table 8** Mean rolling fractal dimensions for the USA (S&P 500) and AUS (S&P/ASX 200) indices across different window lengths

| Market | FD Type | 30-day | 60-day | 100-day |
|---|---|---|---|---|
| USA | Higuchi | 1.5177 | 1.5117 | 1.5076 |
| | Katz | 1.6238 | 1.6251 | 1.6247 |
| | Sevcik | 1.2983 | 1.3098 | 1.3163 |
| AUS | Higuchi | 1.5496 | 1.5358 | 1.5247 |
| | Katz | 1.6411 | 1.6457 | 1.6453 |
| | Sevcik | 1.3147 | 1.3327 | 1.3438 |

While permutation entropy captures changes in the temporal ordering and local unpredictability of price movements, a fuller understanding of market complexity requires examining how the geometric structure of fluctuations evolves across scales. To complement the entropy analysis, we therefore evaluate fractal dimensions using Higuchi's, Katz's, and Sevcik's estimators. These measures provide additional insight into short-, medium-, and long-horizon roughness in the price series, allowing us to assess whether the tariff event altered the broader scaling behaviour of the two markets.

The fractal-dimension results in Table 7 complement the entropy analysis by highlighting how the geometric structure of market fluctuations changed around the tariff event.

The fractal-dimension estimates indicate clear shifts in market complexity following the tariff event. Both indices show an increase in the Higuchi dimension (USA: 1.5072 $\rightarrow$ 1.5612; AUS: 1.5238 $\rightarrow$ 1.5329), reflecting greater short-run irregularity in the post-event period. By contrast, the Katz and Sevcik dimensions decline across both markets, suggesting a smoothing of geometric structure and reduced long-horizon fractal roughness. The adjustment is more pronounced for Australia, particularly in the Sevcik measure (1.3822 $\rightarrow$ 1.2475), indicating a stronger medium-term structural response relative to the USA.

The rolling fractal-dimension estimates in Table 8 following the tariff announcement on 2 April 2025 show clear differences in how the USA and AUS stock markets adjusted to the shock across multiple temporal scales. In the short 30-day window, the AUS market exhibits noticeably higher Higuchi complexity (1.5496) compared with the USA (1.5177), indicating stronger post-event local irregularity and more intense short-horizon fluctuations. A similar pattern emerges in the 60-day window, where the ASX maintains a higher Higuchi dimension (1.5358 vs. 1.5117), and in the 100-day window, where the gap persists (1.5247 vs. 1.5076). These consistently elevated values suggest that the Australian market absorbed the tariff shock through sustained short-run volatility.

**Table 9** Change in fractal dimensions (After − Before) for the USA (S&P 500) and AUS (S&P/ASX 200) indices following the tariff announcement

| Market | FD Type | ΔFD (After − Before) |
|---|---|---|
| USA | Higuchi | +0.0540 |
| | Katz | −0.0475 |
| | Sevcik | −0.1083 |
| AUS | Higuchi | +0.0091 |
| | Katz | −0.0545 |
| | Sevcik | −0.1347 |

The Katz dimension displays an analogous structure: across 30-, 60-, and 100-day windows, the AUS index records higher complexity (1.6411, 1.6457, 1.6453) relative to the USA (1.6238, 1.6251, 1.6247). This indicates that the ASX retained a more geometrically intricate price trajectory during the post-tariff period. By contrast, the Sevcik dimension reveals a more nuanced distinction. For shorter windows, the values for both markets remain relatively close (USA: 1.2983–1.3163; AUS: 1.3147–1.3438), but as the window lengthens, the USA market exhibits a slightly sharper increase in long-run fractal roughness. This suggests that while the AUS market experienced greater short-run complexity following the tariff announcement, the USA developed deeper structural irregularities at longer horizons, reflecting delayed but increasingly pronounced adjustments to the evolving post-tariff environment.

The changes in fractal dimensions presented in Table 9 reveal clear shifts in market complexity following the tariff announcement. Both indices exhibit an increase in the Higuchi dimension (USA: +0.0540; AUS: +0.0091), indicating heightened short-run irregularity, with the response being more pronounced in the USA. In contrast, the Katz and Sevcik dimensions decline across both markets, with especially notable reductions in the Sevcik measure (USA: −0.1083; AUS: −0.1347), suggesting a smoothing of geometric structure and reduced long-horizon fractal roughness. Overall, the USA market shows the stronger rise in short-term complexity, while the Australian market exhibits the more substantial medium-run structural adjustment.

Taken together, the fractal-dimension and permutation-entropy results indicate that the tariff announcement on 2 April 2025 generated only a modest and transient adjustment in the complexity of both markets. The slight increase in the Higuchi dimension suggests a local rise in short-run irregularity, consistent with a brief volatility response, while the decline in Katz and Sevcik dimensions points to a subsequent smoothing of the geometric structure and reduced long-horizon roughness. The permutation-entropy estimates reinforce this view: the small and mostly negative ΔPE values across all windows, paired with near-zero standardized effects, indicate that the underlying unpredictability of ordinal patterns remained broadly stable. Together, these findings imply that although the tariff announcement introduced a mild disturbance, neither the USA nor AUS stock markets experienced a persistent shift in dynamical complexity. Instead, both indices absorbed the shock within a narrow volatility band, returning quickly to pre-event structural behaviour.

## 5 Conclusion, future work and limitation

This study investigates the impact of the April 2, 2025 tariff announcement on the S&P 500 and S&P/ASX 200 indices using a combination of descriptive statistics, machine-learning forecasts, and complexity-based metrics. Market reactions are assessed through changes in permutation entropy and fractal dimensions, which together provide complementary perspectives on short- and medium-term dynamics around the event. Overall, the results point to a muted market response. Entropy measures show only minor variation across all window lengths, indicating that the intrinsic unpredictability of price movements remained largely unchanged. Fractal-dimension estimates reveal similarly modest effects, with a slight increase in short-run irregularity captured by Higuchi's dimension and small declines in Katz's and Sevcik's measures, consistent with temporary fluctuations rather than a structural shift in market behaviour.

Machine-learning results further clarify the nature of price predictability during this period. Ensemble models, particularly Random Forest and XGBoost, consistently outperform linear and distance-based methods, delivering stronger explanatory power and lower forecast errors. In contrast, the weaker performance of kNN and the neural network highlights the importance of methods that balance flexibility and robustness when modelling financial time series subject to event-driven volatility. Taken together, these findings suggest that the tariff announcement did not produce a lasting change in market complexity, while modern ensemble learning techniques remain effective at capturing short-run dynamics even in the presence of policy shocks.

The XAI results indicate that intraday price-level information dominates the predictive structure in both markets, while trading volume and short-horizon directional changes contribute little independent information. The pre–post analysis around April 2, 2025 shows a pronounced increase in the importance of opening prices, especially in the USA, reflecting heightened sensitivity to initial price signals following the event rather than a structural change in price formation dynamics.

Future work could extend the analysis to multiple policy events, sector-level indices, and higher-frequency data to assess whether the observed patterns generalise across different market environments. The study's scope is limited by its reliance on a single event, daily price data, and a relatively small sample, which may restrict the sensitivity of the complexity measures and the full learning capacity of the machine-learning models.

**Author contributions** Conceptualization, Investigation, Methodology, Formal analysis, Writing – original draft, Writing – review & editing:

**Funding** Open Access funding enabled and organized by CAUL and its Member Institutions. Royal Melbourne Institute of Technology.

**Data availability** No datasets were generated or analysed during the current study.

## Declarations

**Conflict of interest** The authors declare that they have no conflict of interest.

**Ethical approval** This article does not contain any studies with human participants or animals performed by any of the authors.